\documentclass[reqno]{amsart}
\usepackage[foot]{amsaddr}

\usepackage[utf8]{inputenc}
\usepackage[T5]{fontenc}

\usepackage{amsmath,amsthm}
\usepackage{amsfonts,amssymb,mathrsfs}

\usepackage[all,cmtip]{xy}

\usepackage[bookmarks=false,breaklinks=true]{hyperref}
\usepackage{float}
\usepackage{graphicx}
\usepackage{color}

\ProvideTextCommand{\DJ}{OT1}{\raisebox{0.25ex}{-}\kern-0.4em D}

\definecolor{myurlcolor}{rgb}{0,0,0.8}
\hypersetup{colorlinks,
linkcolor=myurlcolor,
citecolor=myurlcolor,
urlcolor=myurlcolor}
\usepackage[mathscr]{euscript}

\newcommand{\maps}{\colon}    
\newcommand{\G}{{\mathsf{G}}}   
\renewcommand{\H}{{\mathsf{H}}} 
\newcommand{\M}{{\mathsf{M}}}   
\newcommand{\N}{{\mathsf{N}}}
\newcommand{\Mod}{\mathsf{Mod}}
\newcommand{\Pic}{\mathsf{Pic}} 

\newcommand{\Z}{{\mathbb Z}}  
\newcommand{\R}{{\mathbb R}}  
\newcommand{\C}{{\mathbb C}}  

\newcommand{\To}{\Rightarrow}

\newcommand{\define}[1]{\textbf{\boldmath{#1}}}

\newtheorem{thm}{Theorem}

\theoremstyle{definition}
\newtheorem{defn}[thm]{Definition}

        \newcommand{\be}{\begin{equation}}
        \newcommand{\ee}{\end{equation}}
        \newcommand{\ba}{\begin{eqnarray}}
        \newcommand{\ea}{\end{eqnarray}}
        \newcommand{\ban}{\begin{eqnarray*}}
        \newcommand{\ean}{\end{eqnarray*}}
        \newcommand{\barr}{\begin{array}}
        \newcommand{\earr}{\end{array}}

\begin{document}
\title{Ho\`ang Xu\^an S\'inh's Thesis: \\ Categorifying Group Theory}
\author[Baez]{John C.\ Baez} 
\address{Department of Mathematics, University of California, Riverside CA, 92521, USA}
\address{Centre for Quantum Technologies, National University of Singapore, 117543, Singapore}
\date{\today}
\maketitle

\vbox{\vskip 10 em}
\begin{abstract}
During what Vietnamese call the American War, Alexander Grothendieck spent three weeks teaching mathematics in and near Hanoi.  Ho\`ang Xu\^an S\'inh took notes on his lectures and later did her thesis work with him by correspondence.  In her thesis she developed the theory of `Gr-categories', which are monoidal categories in which all objects and morphisms have inverses.   Now often called `2-groups', these structures allow the study of symmetries that themselves have symmetries. After a brief account of how Ho\`ang Xu\^an S\'inh wrote her thesis, we explain some of its main results, and its context in the history of mathematics.
\end{abstract}

\section{Introduction}

The story of Ho\`ang Xu\^an S\'inh is remarkable because it combines dramatic historical events with revolutionary mathematics.  Some mathematicians make exciting discoveries while living peaceful lives.  Many have their work disrupted or prematurely cut off by wars and revolutions.  But some manage to carry out profound research on the fiery background of  history.

In war-torn Hanoi, Ho\`ang Xu\^an S\'inh met the visionary mathematician Alexander Grothendieck, who had visited to give a series of lectures---in part as a protest against American aggression.  After he returned to France, she did her thesis with him by correspondence, writing it by hand under the light of a kerosene lamp as the bombing of Hanoi reached its peak.  In her thesis she established the most fundamental properties of a novel mathematical structure that takes the concept of symmetry and pushes it to new heights, making precise the concept of \emph{symmetries of symmetries}.   She then traveled to Paris to defend her thesis before some of the most illustrious mathematicians of the era.   After returning to Hanoi, she set up the first private university in Vietnam!   But our story here is about her thesis.

In 1951, Ho\`ang Xu\^an S\'inh completed her bachelor's degree in Hanoi.  She then traveled to Paris for a second baccalaureate in mathematics.  She stayed in France to study for the competitive examination for civil service at the University of Toulouse, which she completed in 1959.  Then she returned to Vietnam and taught mathematics at the Hanoi National University of Education.

In late 1967, when the U.S.\ attacks on Vietnam were escalating, something surprising happened.  Grothendieck visited Vietnam and spent three weeks there teaching mathematics. Ho\`ang Xu\^an S\'inh took the notes for these lectures.  Because of the war, Grothendieck's lectures were held away from Hanoi after the first week---first in the nearby countryside, and later in \DJ\d ai {\fontencoding{T5}T\`\uhorn \selectfont }.   He wrote an account of this which is well worth reading \cite{G67}.  To quote just a bit:
\begin{quote}
Like most more or less public activities, the lectures were scheduled between about
6 and 10 a.m., because the bombing usually took place later in the day, rarely before 11
a.m. During most of my stay the sky was cloud-covered, and consequently there were few
bombing raids. The first serious bombardments had been anticipated; they took place on
Friday 17 November, two days before we left for the countryside. Three times my talk was
interrupted by alarms, during which we took refuge in shelters. Each alert lasted about
ten minutes. Something which is at first very striking to the newcomer is the great calm,
almost indifference, with which the population reacts to the alarms, which have become a
daily routine. 
\end{quote}

At some point Ho\`ang Xu\^an S\'inh asked Grothendieck to be her thesis supervisor, and he accepted.  After he returned to France, she continued to work with him by correspondence.   She later wrote \cite{H23}:

\begin{quote}
If I remember correctly, he wrote to me twice and I wrote to him three times. The first time he wrote to me was to give me the subject of the thesis and the work plan; the second is to tell me to drop the problem of inverting objects if I can't do it. As for me, I think I wrote to him three times: the first was to tell Grothendieck that I couldn't invert objects because of the non-strict commutativity; the second is to tell him that I succeeded in inverting objects; and the third is to tell him that I have finished the job. The letters had to be very short because we were in time of war, eight months for a letter to arrive at its destination between France and Vietnam. When I finished my work in writing, I sent it to my brother, who lives in France, and he brought it to Grothendieck.
\end{quote}

According to the website of Thang Long University, her two main impressions from her contacts with Grothendieck were these:

\begin{enumerate}
\item A good teacher is a teacher who turns something difficult into something easy.

\item We should always avoid anything that is fictitious, live in accordance to our own feelings and value simple people.
\end{enumerate}

She finished her thesis in 1972.   Around Christmas that year, the United States dropped over 20,000 tons of bombs on North Vietnam, mainly Hanoi.  So, it is not surprising that she defended her thesis three years later, when the North had almost won.   But she mentions another reason for the delay \cite{T19}:

\begin{quote}
I was a doctorate student during wartime, bombs and bullets. Back then, I was teaching at Hanoi National University of Education; the school did not provied a way to take leave to study for a doctorate. I taught during the day and worked on my thesis during the night under the kerosene lamp light. I wrote in French under my distant teacher's guidance. When I got the approval from France to come over to defend, there were disagreeable talks about not letting me because they was afraid I wasn't coming back. The most supportive person during the time was Lady H\`a Th\d i Qu{\fontencoding{T5}\selectfont  \'\ecircumflex}---President of the Vietnamese Women Coalescent organization. Madame Qu{\fontencoding{T5}\selectfont  \'\ecircumflex} was a guerilla, without the conditions to get much education, but gave very convincing reasons to support me.  She said, firstly, I was 40 years old, it is very difficult to get a job abroad at 40 years old, and without a job, how can I live? Second, my child is at home, no woman would ever leave her child... so comrades, let's not be worried, let her go.  I finished my thesis in 1972, and three years later with the help and struggle of the women's organization, I was able to travel over to defend in 1975....
\end{quote}

After finishing her thesis, she went to France to defend it at Paris Diderot University (also called Paris VII).  Her thesis committee included not only Grothendieck but also some other excellent mathematicians: 

\begin{itemize}
\item Henri Cartan (a member of Bourbaki famous for his work on sheaf theory, algebraic topology, and potential theory), 
\item Laurent Schwartz (famous for developing the theory of distributions), 
\item Michel Zisman (best known for his work on the ``calculus of fractions'' in category theory), 
\item  Jean-Louis Verdier (a former student of Grothendieck who generalized Poincar\'e duality to algebraic geometry).
\end{itemize}
Her thesis defense lasted two and a half hours.   Soon thereafter she defended a \textsl{second} thesis, entitled \textsl{The Embedding of a One-dimensional Complex in a Two-dimensional Differential Manifold}.  

But it is her first thesis that is our concern here.   Its title was \emph{Gr-cat\'egories}.   But what are Gr-categories, and why are they so important?   I will give a quick answer here, and a more detailed one in the rest of the paper.

Throughout science and mathematics, symmetry plays a powerful simplifying
role.   This is often formalized using group theory.  However, group theory is just the tip of a larger subject that could be called `higher-dimensional group theory'. For example, in many contexts where we might be tempted to use groups, we can use a richer setup where in addition to symmetries of some object, there are higher levels of structure: symmetries of symmetries, symmetries of symmetries of  symmetries, and so on.

For example, if we are studying the rotational symmetries of a sphere, we usually describe a rotation as a special sort of $3 \times 3$ matrix.  But we can also consider a \emph{continuous path} of such rotation matrices, starting at one rotation matrix and ending at another.   Such a path may be considered as a symmetry going between symmetries.

We can continue this pattern to higher levels by considering continuous 1-parameter families of continuous paths of rotation matrices, which describe symmetries of symmetries of symmetries, and so on.    But there is already a fascinating world to discover if we stop at the second level.  This is what the concept of Gr-category formalizes.    Gr-categories got their name because they blend the concepts of \emph{group} and \emph{category}: a Gr-category is a category that resembles a group.  The objects of a Gr-category can be thought of as symmetries, and the morphisms as
symmetries of symmetries.

In Section \ref{sec:monoidal} we describe a key prerequisite for the theory of Gr-categories, namely `monoidal categories', and a bit of the history of this concept, which was still fairly new when Ho\`ang Xu\^an S\'inh wrote her thesis.  In Section \ref{sec:Gr} we explain Gr-categories and Ho\`ang Xu\^an S\'inh's fundamental theorem classifying these.  Then we turn to examples of Gr-categories.   In Section \ref{sec:strict_Gr} we discuss `strict' Gr-categories, where the associative law and unit law hold as equations rather than merely up to isomorphism.  We explain how to get strict Gr-categories from `crossed modules' and use this to give many examples of Gr-categories.  In Section \ref{sec:Pic} we give examples of what could be called `abelian' Gr-categories, which Ho\`ang Xu\^an S\'inh called `Pic-categories' after the mathematician \'Emile Picard.   In Section \ref{sec:topology} we describe how Gr-categories arise in topology, and give more examples using these ideas.   In Section \ref{sec:physics} we give a tiny taste of how Gr-categories---now called `2-groups'---have made an appearance in theoretical and mathematical physics.   Finally, in Appendix \ref{appendix} we give some definitions concerning monoidal and symmetric monoidal categories, allowing us to state some of Ho\`ang Xu\^an S\'inh's results more precisely.

\section{Monoidal categories}
\label{sec:monoidal}

Ho\`ang Xu\^an S\'inh's work on Gr-categories is part of the trend toward `categorification': taking familiar structures built using sets and finding their analogues built using categories.   As we shall see, a  Gr-category is a category that resembles a group.  But a crucial subtlety is that a Gr-category has inverses at two different levels: both the objects and morphisms have inverses!  The objects describe symmetries, while the morphisms describe symmetries \emph{of} symmetries.

To define Gr-categories we should start with monoidal categories.  A monoidal category is a category that resembles a monoid.   A bit more precisely, a monoidal category is a category $\M$ with a tensor product
\[           \otimes \colon \M \times \M \to \M \]
and a unit object $I \in \M$, obeying  the associative and unit laws up to specified isomorphisms, which in turn must obey some laws of their own.   Given this, a Gr-category is a monoidal category in which every morphism and every object has an inverse.   Here a morphism $f$ has an inverse if there is a morphism $f'$ such that $f \circ f'$ and $f' \circ f$ are identity morphisms.  Similarly, an object $g$ in a monoidal category has an inverse if there an object $g'$ such that $g \otimes g' \cong I$ and $g' \otimes g \cong I$.

Now that monoidal categories are well-understood, the definition of Gr-category seems simple and natural.   However, monoidal categories were still quite new when Ho\`ang Xu\^an S\'inh wrote her thesis, and she did even not call them that.  In 1963, B\'enabou published a short paper containing a preliminary version of monoidal categories under the name of `cat\'egories avec multiplication' \cite{B63}.   Also in 1963, Mac Lane published a paper correcting an important problem with B\'enabou's definition and giving the modern definition  \cite{M63}.  However, at this time Mac Lane called monoidal categories `bicategories'.  Ironically, this term is now used for a very different concept introduced later by B\'enabou \cite{B67}.   According to Mac Lane \cite{M13}, it was Eilenberg who coined the term `monoidal category'.  It seems to have first appeared in a paper by Eilenberg and Kelly \cite{EK}.   Ho\`ang Xu\^an S\'inh does not cite Mac Lane's 1963 paper.  For her terminology on monoidal categories, she refers to the thesis of Neantro Saavedra-Rivano, another student of Grothendieck \cite{SR70,SR72}.

But let us see how B\'enabou handled the concept of `category with multiplication', and how Mac Lane's important correction plays a key role in Ho\`ang Xu\^an S\'inh's work.   First, B\'enabou demanded that a category with multiplication have a functor
\[           \otimes \colon \M \times \M \to \M. \]
We call this functor the \define{tensor product}, and write $\otimes(x,y)=x \otimes y$ and $\otimes(f,g)=f \otimes g$
for objects $x, y \in M$ and morphisms $f, g$ in $M$.  Second, B\'enabou said that a category with multiplication have an object $I \in \M$ called the \define{unit}.  For example, the unit for the tensor product in the category of vector spaces can be taken to be any one-dimensional vector space, while the unit for the Cartesian product in the category of sets can be taken to be any one-element set.  

B\'enabou did not demand that the tensor product be
associative `on the nose', in the obvious equational way:
\[         (x \otimes y) \otimes z = x \otimes (y \otimes z) .\]
Instead, he required that the tensor product be associative up to a natural isomorphism, which we now call the `associator':
\[     a_{x,y,z} \maps (x \otimes y) \otimes z \xrightarrow{\;\sim\;} x \otimes (y \otimes z)  .\]
Similarly, he did not demand that $I$ act as the unit for
the tensor product on the nose, but only up to natural isomorphisms
which are now called the left and right `unitors':
\[       \ell_x \maps I \otimes x \xrightarrow{\;\sim\;} x  , \qquad
            r_x \maps x \otimes I \xrightarrow{\;\sim\;} x  .\]
The reason is that in applications, it is usually too much to
expect equations between objects in a category: usually we just
have isomorphisms, and this is good enough!  Indeed this is a
basic principle of categorification: equations between objects are
bad; we should instead specify isomorphisms.

If we stopped the definition here, there would be a problem, since we could
use the associator to build \emph{several different} isomorphisms
between a pair of objects constructed using tensor products.  The simplest
example occurs when we have four objects $w,x,y,z \in \M$.  Then we can
build two isomorphisms from $((w \otimes x) \otimes y) \otimes z$ to $w \otimes (x \otimes (y \otimes z))$ using associators:
\be
\xy 0;/r.3pc/:
    (-24.73,8.03)*+{ (w \otimes (x \otimes y)) \otimes z}="l";
    (0,26)*+{ ((w \otimes x) \otimes y) \otimes z}="t";
    (24.73,8.03)*+{ (w \otimes x) \otimes (y \otimes z)}="r";
    (15.28,-21.03)*+{ w \otimes (x \otimes (y \otimes z))}="br";
    (-15.28,-21.03)*+{ w \otimes ((x \otimes y) \otimes z)}="bl";
     {\ar_{a_{w, x, y} \otimes 1_z} "t";"l"};
     {\ar_{a_{w, x \otimes y, z}} "l";"bl"};
     {\ar_{1_w \otimes a_{x,y,z}} "bl";"br"};
     {\ar^{a_{w, x, y \otimes z}} "r";"br"};
     {\ar^{a_{w \otimes x, y, z}} "t";"r"};
\endxy
\label{pentagon}
\ee
If these isomorphisms were different, we would need to to say \emph{how} we reparenthesized one expression to get the other, just to specify an isomorphism between them.  This situation is theoretically possible, but doesn't seem to come up much in mathematics.   So, it seems wise to require that the above pentagon commutes.  It does in most examples we usually care about.

When we build isomorphisms using unitors, we get further
ways to build multiple isomorphisms between objects constructed using tensor products and the unit object.  The simplest example is this:
\be 
\xymatrix{
I \otimes I   \ar@/^/[rr]^{\ell_I}  \ar@/_/[rr]_{r_I} && I
}  
\label{bigon}
\ee
When we combine unitors and associators, the three simplest examples are these triangles:
\be
\xy 0;/r.3pc/:
(-14,10)*+{(I \otimes x) \otimes y}="l";
(14,10)*+{ I \otimes (x \otimes y)}="r";
(0,0)*+{x \otimes y }="b";
 {\ar^{a_{I,x,y}} "l";"r"};
{\ar_{\ell_x \otimes 1_y} "l";"b"};
{\ar^{\ell_{x \otimes y}} "r";"b"};
\endxy    
\label{triangle1}
\ee
\be
\xy 0;/r.3pc/:
(-14,10)*+{(x \otimes I) \otimes y}="l";
(14,10)*+{x \otimes (I \otimes y)}="r";
(0,0)*+{x \otimes y }="b";
 {\ar^{a_{x,I,y}} "l";"r"};
{\ar_{r_x \otimes 1_y} "l";"b"};
{\ar^{1_x \otimes \ell_y} "r";"b"};
\endxy    
\label{triangle2}
\ee
\be
\xy 0;/r.3pc/:
(-14,10)*+{(x \otimes y) \otimes I}="l";
(14,10)*+{x \otimes (y \otimes I)}="r";
(0,0)*+{x \otimes y }="b";
 {\ar^{a_{x,y,I}} "l";"r"};
{\ar_{r_{x \otimes y}} "l";"b"};
{\ar^{1_x \otimes r_y} "r";"b"};
\endxy    
\label{triangle3}
\ee
So, we should require that these triangles commute too.

Of course, there are infinitely many other diagrams we can build using associators and unitors.   B\'enabou tried to handle all of them in one blow.   He imposed an axiom saying essentially that \emph{all} diagrams constructed from associators, left unitors and right unitors must commute.   Unfortunately such an axiom fails to have the desired effect unless one states it very carefully---which B\'enabou, alas, did not.  The problem is that even in a perfectly nice monoidal category, there can be diagrams built from the associator that fail to commute due to `accidental' equations between tensor products.  For instance, suppose in some monoidal category we just happen to have $(x\otimes y)\otimes z = b \otimes (c \otimes d) $ and $(b \otimes c) \otimes d = (u \otimes v) \otimes w$ and $x\otimes (y\otimes z) = u \otimes (v\otimes w)$.  Then we can build this diagram:
\[
\xy 0;/r.3pc/:
(-24,15)*+{(u \otimes v) \otimes w = (b \otimes c) \otimes d}="l";
(24,15)*+{b \otimes (c \otimes d) = (x \otimes y) \otimes z}="r";
(0,0)*+{u \otimes (v \otimes w) = x \otimes (y \otimes z)}="b";
 {\ar^{a_{b,c,d}} "l";"r"};
{\ar_{a_{u,v,w}} "l";"b"};
{\ar^{a_{x,y,z}} "r";"b"};
\endxy    
\]
According to B\'enabou's axiom, this diagram must commute.  But this does not hold in many of the examples we are interested in.   So, B\'enabou's axiom is too strong.

Mac Lane's solution \cite{M63} was to show that if diagrams \eqref{pentagon}, \eqref{bigon}, \eqref{triangle1}, \eqref{triangle2} and \eqref{triangle3} commute, all diagrams built from associators and unitors that we really \emph{want} to commute actually \emph{do} commute.  Making this precise and proving it was a major feat: it is called Mac Lane's coherence theorem.   In 1964, Kelly \cite{K64} improved Mac Lane's result by showing that the pentagon identity \eqref{pentagon} and the middle triangle identity \eqref{triangle2} imply the rest.   Thus, the modern definition of monoidal category is as follows:

\begin{defn}
 A \textbf{monoidal category} is a category $M$ together with a \define{tensor
product} functor $\otimes \maps M \times M \to M$, a \define{unit object}
$I \in M$, and natural isomorphisms
called the \textbf{associator}:
\[ a_{x,y,z} \maps (x \otimes y) \otimes z \xrightarrow{\; \sim \;} x \otimes (y \otimes z), \]
the \textbf{left unitor}:
\[ \ell_x \maps I \otimes x \xrightarrow{\; \sim \;} x , \]
and the \textbf{right unitor}:
\[ r_x \maps x \otimes I \xrightarrow{\; \sim \;} x. \]
such that the pentagon \eqref{pentagon} and the triangle \eqref{triangle2} commute.
\end{defn}

As noted by Mac Lane, the pentagon was also discussed in a 1963 paper by James Stasheff \cite{S63}, as part of an infinite sequence of polytopes called
`associahedra'.  Stasheff defined a concept of `$A_\infty$-space',
which is roughly a topological space having a product that is
associative up to homotopy, where this homotopy satisfies the pentagon
identity up homotopy, that homotopy satisfies yet another identity up
to homotopy, and so on, \textit{ad infinitum}.  The $n$th of these
identities is described by the $n$-dimensional associahedron.  The
first identity is just the associative law, which plays a crucial role
in the definition of monoid.  Mac Lane realized that the second, the pentagon
identity, should play a similar role in the definition of monoidal
category.  The higher ones show up, at least implicitly, in Mac Lane's
proof of his coherence theorem.  They become even more important when
we consider monoidal bicategories, monoidal tricategories and so on.  

Less glamorous than the associator, the unitors and the laws they obey
are also crucial to Mac Lane's coherence theorem.  Together with the associahedra, they are described by a sequence of cell complexes sometimes called `monoidahedra' \cite{T99} or `unital associahedra' \cite{MT11}.    

\section{Gr-categories}
\label{sec:Gr}

We have already stated Ho\`ang Xu\^an S\'inh's definition of Gr-category, but now we are ready to examine her work on this concept more carefully, so let us repeat the definition:

\begin{defn}
A \define{Gr-category} is a monoidal category such that every object has an inverse and every morphism has an inverse.
\end{defn}

She not only defined Gr-categories; she also defined when two Gr-categories are `equivalent'.    The idea is that equivalent Gr-categories are the same for all practical purposes.   The definition of equivalence is a bit lengthy, so we relegate it to Appendix \ref{appendix}.   But the idea, even in a rough form, is very interesting.  

We can define maps between Gr-categories.   We say two Gr-categories $\G$ and $\G'$ are isomorphic if there are maps $F \maps \G \to G'$ and $G \maps \G' \to \G$ that are inverses, so
\[    GF = 1_\G \textrm{ \; and \; } FG = 1_{\G'}   .\]
This is of course a typical definition of isomorphism for any sort of mathematical gadget.  But because Gr-categories are categories, it turns out we can also define isomorphisms \emph{between} maps between Gr-categories.  We then say two Gr-categories $\G$ and $\G'$ are `equivalent' if there are maps $F \maps \G \to G'$ and $G \maps \G' \to \G$ that are inverses up to isomorphism:
\[    GF \cong 1_\G \textrm{ \; and \; } FG \cong 1_{\G'}   .\]

Two groups are isomorphic if renaming the elements of one group gives the other.  Similarly, two Gr-categories are isomorphic if and only if we can get from one to the other by renaming the objects and morphisms.  But equivalence is more flexible!   For example, suppose we have a Gr-category $\G$ with a collection of isomorphic objects.  Then we can \emph{remove} all but one of these and get a smaller Gr-category that is equivalent to $\G$.  The idea is that it counts as `redundant' to have multiple isomorphic copies of an object in a Gr-category.  We can remove these redundant copies and get an equivalent Gr-category.

It is worth spelling this out in a bit more detail in an example.  Suppose a Gr-category $\G$ has an isomorphism between two distinct objects, say $f \maps g \xrightarrow{\, \sim \, } g'$.  Then we can remove $g$ from $\G$, along with all the isomorphisms going to and from it, and get a smaller category $\G'$.  This category is not immediately a Gr-category, since the tensor product of two objects in $\G'$ may happen to equal $g$.  But we can redefine such tensor products to equal $g'$, and suitably `reroute' the associator and unitors of $\G'$ using the isomorphism $f$, to make $\G'$ into a Gr-category.   Furthermore, this new Gr-category $\G'$ will be equivalent to $\G$.

More generally, we can show any Gr-category is equivalent to a \define{skeletal} Gr-category: one where any two isomorphic objects must be equal.   To classify Gr-categories up to equivalence, it thus suffices to classify skeletal Gr-categories up to equivalence.

Working with skeletal Gr-categories simplifies the classification problem.   In a general Gr-category, the group laws hold up to isomorphism, since
\[    (g \otimes g') \otimes g'' \cong g \otimes (g' \otimes g'') \]
\[      I \otimes g \cong g \cong g \otimes I  \]
and for any object $g$ there is an object $g'$ with
\[      g \otimes g' \cong I \cong g' \otimes g. \]
But in a skeletal Gr-category, isomorphic objects are equal, so these laws become equations.  Thus in a skeletal Gr-category the set of objects forms a group!  

Indeed, from a skeletal Gr-category $\G$ we easily get three pieces of data:
\begin{itemize}
\item the group $G$ of objects,
\item the group $A$ of automorphisms of the unit object $I \in \G$, which is an abelian group,
\item an action $\rho$ of the group $G$ as automorphisms of the abelian 
group $A$, given as follows:
\[     \rho(g)a := 1_g \otimes a \otimes 1_{g^{-1}} .\]
\end{itemize}

Clearly, the group $G$ tells everything about tensoring objects in $\G$.  More subtly, the group $A$ tells us everything we need to know about composing morphisms in $\G$.  Why is this?  First, in a skeletal Gr-category every morphism is an isomorphism, and isomorphic objects are equal.  Thus every morphism in $\G$ is an automorphism $f \maps g \to g$ for some $g \in \G$.   Second, the group of automorphisms of any $g \in \G$ is 
isomorphic to the group of automorphisms of $I \in \G$, via the map sending $f \maps g \to g$ to $1_{g^{-1}} \otimes f \maps I \to I$.  

The group $A$ is abelian because given $a, b \maps I \to I$ we have
\[   ab =  (a 1_I) \otimes (1_I b) = (a \otimes 1_I)(1_I \otimes b) = a \otimes b = (1_I \otimes a)(b \otimes 1_I) = (1_I b) \otimes (a 1_I) = ba   .\]
This marvelous calculation, later called the `Eckmann--Hilton argument' \cite{EH62},
also shows that tensoring automorphisms of $I \in \G$ is the same as composing them!  A further argument shows that the rest of the information about tensoring morphisms in $\G$ is contained in the action $\rho$ of $G$ on $A$.

The elephant in the room is the associator.  To describe the associator of a Gr-category
we need a fourth piece of data---the subtlest and most interesting piece.   As we have seen, in any monoidal category the tensor product is associative up to a natural transformation called the associator:
\[    \alpha_{g_1, g_2, g_3} \maps (g_1 \otimes g_2) \otimes g_3 \to g_1 \otimes (g_2 \otimes g_3).  \]
Since this takes three objects and gives a morphism, we should try to encode this into a
map from $G^3$ to $A$.  How can we do this?

In a skeletal Gr-category, isomorphic objects are equal, so we can drop the parentheses and write
\[    \alpha_{g_1, g_2, g_3} \maps g_1 \otimes g_2 \otimes g_3 \to g_1 \otimes g_2 \otimes g_3 .\]
But beware: $\alpha_{g_1,g_2,g_3}$ may not be an identity morphism!  It is just some automorphism of $ g_1 \otimes g_2 \otimes g_3$.  We can turn this into an automorphism of the unit object $I$ and define a map
\[  a \maps G^3 \to A \]
as follows:
\[    a(g_1, g_2, g_3) = 1_{(g_1 \otimes g_2 \otimes g_3)^{-1}} \otimes \alpha_{g_1, g_2, g_3} .\]

Not every map $a \maps G^3 \to A$ arises from some Gr-category.  After all, the associator must obey the pentagon identity \eqref{pentagon}.   In terms of the map $a$, this says that for all $g_1, g_2, g_3, g_4 \in G$ we have
\be
  \rho(g_0)a(g_1,g_2,g_3) - a(g_0g_1, g_2, g_3) +  
  a(g_0, g_1g_2, g_3) - a(g_0, g_1, g_2g_3) +
  a(g_0, g_1, g_2) = 0.
\label{3-cocycle}
\ee
Here we have switched to writing the group operation in $A$ additively, since this group is abelian.

Now, the precise value of the map $a \maps G^3 \to A$ depends on our choice of a skeletal Gr-category equivalent to $\G$.   Thus $a$ is not uniquely determined by $\G$.  With some calculation, we can check that by changing our choice we can change the map $a\maps G^3 \to A$ to any map $a' \maps G^3 \to A$ with 
\[          a' - a = df   \]
where $f \maps G^2 \to A$ is an arbitrary map and
\be     (df)(g_1, g_2, g_3) = \rho(g_1)f(g_2, g_3) - f(g_1 g_2, g_3) + f(g_1, g_2 g_3) - f(g_1, g_2) . 
\label{3-coboundary}
\ee

Equations  \eqref{3-cocycle} and \eqref{3-coboundary} look intimidating at first sight, but by the time Ho\`ang Xu\^an S\'inh wrote her thesis these equations were already familiar in the subject of group cohomology. For any group $G$ equipped with an action on an abelian group $A$, there is a sequence of groups $H^n(G,A)$ called the `cohomology groups' of $G$ with coefficients in $A$ \cite{R79}.   For $n = 1$ and $n = 2$ these groups first arose in Galois theory and the study of group extensions, and they are fairly easy to interpret.  For $n = 3$ and above, these groups were independently discovered in 1943--45 by teams in different countries who found it difficult to communicate during the chaos of World War II.  For the complex history of these developments see the accounts of Hilton \cite{H02}, Mac Lane \cite{M78,M79,M88} and Weibel \cite{W99}.

The significance of these higher cohomology groups, even $H^3(G,A)$, was at first rather obscure.  Around 1945, Eilenberg and Mac Lane \cite{EM45} found a powerful topological explanation of these groups.  Ultimately, the idea is that that for any group $G$ there is a topological space called the `Eilenberg--Mac Lane space' $K(G,1)$ whose cohomology groups with a suitable system of `local coefficients' depending on $A$ and the action $\rho$ are the groups $H^n(G,A)$.  But the meaning of $H^3(G,A)$ became even more vivid thanks to Ho\`ang Xu\^an S\'inh's work.

How does this come about?  First, the resemblance of equations \eqref{3-cocycle} and \eqref{3-coboundary} is no coincidence!  For any map $h \maps G^n \to A$ there is a similar-looking formula for a map $dh \maps G^{n+1} \to A$.  A map $h \maps G^n \to A$ with $dh = 0$ is called an \define{$n$-cocycle}, and a map $h \maps G^n \to A$ of the form $h = df$ for some $f$ is called an \define{$n$-coboundary}.   Since 
\[           d d f = 0 \]
for all $f$, every $n$-boundary is automatically an $n$-cocycle.    The group of $n$-cocycles modulo $n$-coboundaries is called the \define{$n$th cohomology group} of $G$ with coefficients in $A$, and denoted $H^n(G,A)$.   Elements of this group are called \define{cohomology classes}.

In this language, \eqref{3-cocycle} says that the associator of a skeletal Gr-category gives a 3-cocycle $a \maps G^3 \to A$.   Similarly, \eqref{3-coboundary} says that the 3-cocycles $a$ and $a'$ differ by a 3-coboundary---and thus represent the same cohomology class.  

Thus, Ho\`ang Xu\^an S\'inh's work shows that any Gr-category gives, not only two groups $G$ and $A$ and an action of the first on the second, but also an element of $H^3(G,A)$. More importantly, she proved that two Gr-categories with the same groups $G$ and $A$ and the same action $\rho$ of $G$ on $A$ are equivalent \emph{if and only if} their associators define the same element of $H^3(G,A)$.  We state this result more carefully in Theorem \ref{thm:Gr1} of Appendix \ref{appendix}. 

This gave a new explanation of the meaning of the cohomology group $H^3(G,A)$.  In simple terms, this group classifies the possible associators that a Gr-category can have when the rest of its structure is held fixed.  The element of $H^3(G,A)$ determined by the associator of a Gr-category $\G$ is now called the \define{S\'inh invariant} of $\G$.

\section{Strict Gr-categories}
\label{sec:strict_Gr}

We have emphasized that in a Gr-category, the associative and unit laws hold only up to natural isomorphism.   Nonetheless it is interesting to consider Gr-categories where these natural isomorphisms are identity morphisms.   In a skeletal Gr-category we have
\[    (g \otimes g') \otimes g'' = g \otimes (g' \otimes g'') \]
and
\[      I \otimes g = g = g \otimes I.  \]
However, this is somewhat deceptive, because the associator and unitors may still not be identity morphisms!   Thus we say a Gr-category is \define{strict} if its associator and unitors are identity morphisms.   

In 1978, Ho\`ang Xu\^an S\'inh published a paper about strict Gr-categories in which she proved that every Gr-category is equivalent to a strict one \cite{H78}.  At a first glance this might seem to contradict that fact that every Gr-category $\G$ is equivalent to a skeletal one, which is then characterized up to equivalence by four pieces of data:
\begin{itemize}
\item its group $G$ of objects,
\item the abelian group $A$ of automorphisms of the unit object,
\item the action $\rho$ of $G$ on $A$,
\item the cohomology class $[a] \in H^3(G,A)$ arising from its associator.
\end{itemize}
Since any Gr-category is equivalent to a strict one, where the associator is the identity, isn't $[a] \in H^3(G,A)$ always zero?  No, because that strict Gr-category is not generally skeletal!   In fact, it follows from what we have said that $[a] = 0$ if and only if $\G$ is equivalent to a strict and skeletal Gr-category.  

This fact bears repeating, since it trips up beginners so often: \emph{every Gr-category is equivalent to a strict one, and also equivalent to a skeletal one, but it is equivalent to a strict and skeletal one if and only if its S\'inh invariant $[a]$ is zero}.

Before Gr-categories were introduced, strict Gr-categories were already known, though of course not under that name.  Their early history is quite obscure, at least to this author. Strict Gr-categories are now often called `categorical groups', but the earliest paper I have found using that term dates only to 1981, and uses it to mean something else \cite{S81}.   In 1976, Brown and Spencer \cite{BS76} published a proof that strict Gr-categories, which they called `$\mathcal{G}$-groupoids', are equivalent to `crossed modules', a structure introduced by J.\ H.\ C.\ Whitehead in 1946 for use in topology \cite{W46,W49}.   But they write that this result was known to Verdier in 1965, used by Duskin in some unpublished notes in 1969, and discovered independently by them in 1972.  

Crossed modules are still a useful way to get examples of Gr-categories. For a modern proof that crossed modules are equivalent to strict Gr-categories, see Forrester-Barker \cite{FB}.  But what exactly is a crossed module?

\begin{defn}
A \define{crossed module} is a quadruple $(G,H,t,\rho)$ where $G$ and $H$ are groups, $t \maps H \to G$ is a homomorphism, and $\rho$ is an action of $G$ as automorphisms of $H$ such that $t$ is \define{$G$-equivariant}:
\[   t(\rho(g)h) = g \, t(h)\, g^{-1} \]
and $t$ satisfies the so-called \define{Peiffer identity}:
\[    \rho(t(h)) h' = hh'h^{-1} . \]
\end{defn}

We can obtain a crossed module from a strict Gr-category $\G$ as follows:
\begin{itemize}
\item Let $G$ be the set of objects of $\G$, made into a group by tensor product.
\item Let $H$ be the set of all morphisms from $I \in \G$ to arbitrary objects of $\G$, made into a group by tensor product.
\item Let $t \maps H \to G$ map any morphism $f \maps 1 \to g$ to $g \in \G$.
\item Let $\rho(g)h = 1_g \otimes h \otimes 1_{g^{-1}}$.
\end{itemize}

Conversely, we can build a strict Gr-category $\G$ from a crossed module $(G,H,t,\rho)$ as follows.   We take the set of objects of $\G$ to be $G$, and we define the tensor product of objects using multiplication in $G$.  We take the set of morphisms of $\G$ to be the semidirect product $H \rtimes G$ in which multiplication is given by
\[    (h,g) (h',g') = (h \rho(g)h', gg') .\]
We define the tensor product of morphisms using multiplication in this semidirect product.   We treat the pair $(h,g)$ as a morphism from $g \in \G$ to $t(h) g \in \G$, and we define the composite of morphisms
\[  (h,g) \maps g  \to g'  , \qquad \quad
(h',g') \maps g' \to g'' \] 
to be
\[   (hh',g) \maps g \to g'' .\]

Using this correspondence we can now get many examples of strict Gr-categories.   For example: any action of a group on an abelian group, any normal subgroup of a group, or any central extension of a group gives a crossed module, and thus a strict Gr-category.   So, strict Gr-categories unify a number of important concepts in group theory!

To see this, first suppose $\rho$ is any action of a group $G$ on an abelian group $H$. Then defining $t \maps H \to G$ by $t(h) = 1$, the $G$-equivariance of $t$ and the Peiffer identity are easy to check, so we get a crossed module.  Because $t(h) = 1$, the resulting strict Gr-category is skeletal.  In fact we can turn this argument around, and show that any skeletal strict Gr-category arises this way (up to isomorphism).   

Second, suppose $H$ is any normal subgroup of $G$.  Let
$t \maps H \to G$ be the inclusion of $H$ in $G$ and let 
\[   \rho(g) h = g h g^{-1}  \]
Then clearly $t$ is $G$-equivariant and the Peiffer identity holds, so
we get a crossed module.  Here the resulting strict Gr-category is usually
\emph{not} skeletal, since typically we do not have $t(h) = 1$ for all
$h \in H$.

Third, suppose we have a central extension of a group $G$ by a group
$C$, i.e.\ a short exact sequence
\[       1 \xrightarrow{\; i \; } C \xrightarrow{\phantom{\; i \; }} H \xrightarrow{\; t \; } G \xrightarrow{\phantom{\; t \; }} 1 \]
where the image of $C$ lies in the center of $H$.   Then we can define
an action $\rho$ of $G$ on $H$ by choosing any function $j \maps
G \to H$ with $t (j(g)) = g$ for all $g \in G$ and letting
\[            \rho(g) h = j(g) h j(g)^{-1}   .\]
We can check that $\rho$ does not depend on the choice of $j$.   
With more work, we can check that $(G,H,t,\rho)$ is a crossed module.
The resulting strict Gr-category is again usually not skeletal.  On the contrary, since $t$ is surjective, in this Gr-category all objects are isomorphic.

If we build a strict Gr-category from a crossed module, and it is not skeletal,
it may have a nontrivial S\'inh invariant.    For the story of how mathematicians discovered the relation of the 3rd cohomology of groups to crossed modules, see Mac Lane's historical remarks \cite{M79,M88}.  The connection was certainly known well before Ho\`ang Xu\^an S\'inh wrote her thesis.

\section{Pic-categories}
\label{sec:Pic}

The subject of Gr-categories is interesting because there are examples arising from many different branches of mathematics, and the relationships between these examples exposes links between these branches.  One important and tractable class of Gr-categories studied by Ho\`ang Xu\^an S\'inh are what she called `Pic-categories'.  Let us take a look at these.

Any monoidal category $\M$ gives a Gr-category called its \define{core}, namely the subcategory consisting of only the invertible objects and only the invertible morphisms.   For example, if $R$ is any commutative ring, there is a category $R \Mod$ where
\begin{itemize}
\item objects are $R$-modules,
\item morphisms are homomorphisms of $R$-modules.
\end{itemize}
This category becomes monoidal with the usual tensor product of $R$-modules, 
and the unit for this tensor product is $R$ itself, considered
as an $R$-module.  The core of $R \Mod$ is a Gr-category that could
be called $\Pic(R)$, the \define{Picard category} of $R$.  By Ho\`ang Xu\^an S\'inh's
classification of Gr-categories, all the information in this Gr-category is contained in four pieces of
data:
\begin{itemize}
\item 
The group of isomorphism classes of invertible $R$-modules.  This is
called the \define{Picard group} $\mathrm{Pic}(R)$.  This group is always abelian.
\item 
The abelian group of invertible $R$-module homomorphsms $f \maps R \to R$.
Since these are all given by multiplication by invertible elements of $R$, this is
called the \define{units group} $R^\times$.
\item
The action of the Picard group on the units group.  This action is always trivial.
\item
The S\'inh invariant of $\Pic(R)$.  This is always trivial.
\end{itemize}
It is evident from this that the Picard category of a commutative ring is a specially simple sort of Gr-category.

The Picard group and units group are important invariants of a commutative
ring.  For example, when $R$ is the ring of algebraic integers of some algebraic number field, the Picard group $\mathrm{Pic}(R)$ is isomorphic to the `ideal class group' of $R$, familiar in number theory.  The reason is that in this case, every invertible $R$-module is isomorphic to an ideal of $R$.  Even better, it turns  out that for a ring of algebraic integers $\mathrm{Pic}(R)$ is trivial if and only if unique prime factorization holds in $R$---at least up to reordering and units.   When $\mathrm{Pic}(R)$ is nontrivial, unique factorization fails, and this is one reason Dedekind, building on earlier work of Kummer, introduced the ideal class group.  We can compare two examples:

\begin{itemize}
\item
When $R = \Z[i]$ is the ring of Gaussian integers, the Picard group 
$\Pic(R)$ is trivial and the units group $\R^\times$ is $\Z/4$, consisting of
$1, i, -1$ and $-i$.  We have unique prime factorization in $\Z[i]$.
\item
When $R = \Z[\sqrt{-5}]$, the Picard group $\Pic(R)$ is $\Z/2$ and the
units group $R^\times$ is $\Z/2$, consisting of $1$ and $-1$.  We do not have
unique prime factorization in $\Z[\sqrt{-5}]$.
\end{itemize}

It is worth comparing an example from topology.  Suppose $X$ is a compact Hausdorff space.  Then we can let $R$ be the ring of continuous complex-valued functions on $X$, with pointwise addition and multiplication.  In this case $\Pic(R)$ is equivalent to the Gr-category where
\begin{itemize}
\item objects are complex line bundles over $X$,
\item morphisms are isomorphisms between complex line bundles,
\end{itemize}
and the tensor product is the usual tensor product of line bundles.
It follows that we can identify the Picard group $\mathrm{Pic}(R)$ with the set of
isomorphism classes of complex line bundles over $X$, made into a group using
tensor products.   This is famously isomorphic to the cohomology group $H^2(X,\Z)$.  The units group $R^\times$ is simply the group of continuous functions from $X$ to $\C^\times$, the multiplicative group of invertible complex numbers. 

Line bundles also show up in another example of a Gr-category.  
Suppose $X$ is a complex projective algebraic variety.
Then the category of holomorphic vector bundles over $X$, with its usual tensor
product, is a monoidal category.  Its core is a Gr-category where
\begin{itemize}
\item objects are holomorphic line bundles over $X$,
\item morphisms are isomorphisms between holomorphic line bundles,
\end{itemize}
and the tensor product is the usual tensor product of holomorphic line bundles. In this case the group of isomorphism classes of objects is called the \define{Picard group} of $X$.   This Picard group is much more interesting than the previously mentioned purely topological example.  It depends not only on the topology of $X$, but on its holomorphic structure.  Moreover, instead of a discrete group, this Picard group is best thought of as a topological group whose connected components are themselves projective algebraic varieties!

All the examples of Gr-categories mentioned in this section so far are not only monoidal categories, but `symmetric' monoidal categories.  These were introduced by Mac Lane in his 1963 paper \cite{M63}, though not under that name.  Just as monoidal categories are a categorification of monoids, symmetric monoidal categories are a categorification of \emph{commutative} monoids.   But instead of requiring that the tensor product commute `on the nose', we demand that it commute up to a natural isomorphism, which must obey some laws of its own:

\begin{defn}  A \define{symmetric monoidal category} is a monoidal category $\M$ equipped with a natural isomorphism called the \define{symmetry}
\[        S_{x,y} \maps x \otimes y \to y \otimes x  \]
such that 
\[S_{y,x} \circ S_{x,y} = 1_{x \otimes y} \]
for all objects $x,y \in \M$ and the following diagram commutes for all $x,y,z \in \M$:
\be
\xy
   (-15,26)*+{(x \otimes  y) \otimes  z}="tl";
   (15,26)*+{(y \otimes  x) \otimes  z}="tr";
   (-30,0)*+{x \otimes  (y \otimes  z)}="ml";
   (30,0)*+{y \otimes  (x \otimes  z)}="mr";
   (-15,-26)*+{(y \otimes  z) \otimes  x}="bl";
   (15,-26)*+{y \otimes  (z \otimes  x)}="br";
        {\ar^{ S_{x,y} \otimes  z} "tl";"tr"};
        {\ar^{ a_{y,x,z}} "tr";"mr"};
        {\ar^{ a ^{-1}_{x,y,z}} "ml";"tl"};
        {\ar_{ S_{x,y \otimes  z}} "ml";"bl"};
        {\ar_{ a_{y,z,x}} "bl";"br"};
        {\ar^{ y \otimes   S_{x,z} } "mr";"br"};
\endxy
\label{hexagon} 
\ee
\end{defn}
This commuting hexagon says that switching the object $x$ past
$y \otimes z$ all at once is the same as switching it first past $y$ and then
past $z$ (with some associators thrown in to move the parentheses).

Ho\`ang Xu\^an S\'inh made the following definition:
\begin{defn}
A \define{Pic-category} is a symmetric monoidal Gr-category.
\end{defn}
In fact, all the examples of Gr-categories we have seen so far are even better than Pic-categories: they are what Ho\`ang Xu\^an S\'inh called `restreintes'.

\begin{defn}
A Pic-category is \define{restrained} if the symmetry $S_{g,g} \maps g \otimes g \to g \otimes g$ is the identity for every object $g \in \G$.
\end{defn}
\noindent
In particular, what we are calling the Picard category of a commutative ring and the Picard category of a complex projective algebraic variety are restrained Pic-categories.

Among Gr-categories, restrained Pic-categories are especially simple.  Recall that Gr-categories are classified up to equivalence by four pieces of data:
\begin{itemize}
\item the group $G$ of isomorphism classes of objects,
\item the abelian group $A$ of automorphisms of the unit object,
\item an action $\rho$ of $G$ as automorphisms of $A$,
\item an element $[a] \in H^3(G,A)$.
\end{itemize}
A restrained Pic-category always has these properties:
\begin{itemize}
\item the group $G$ is abelian,
\item the action $\rho$ is trivial,
\item the element $[a]$ is zero.    
\end{itemize}
In fact, Ho\`ang Xu\^an S\'inh showed in her thesis that any restrained Pic-category is characterized by a pair of abelian groups $G$ and $A$.    To do this, she first defined an appropriate concept of `equivalence' for Pic-categories, which is more fine-grained than equivalence of Gr-categories since the symmetry also plays a role.  Then, she classified Pic-categories and also restrained Pic-categories up to this notion of equivalence.   We explain the latter classification in Theorem \ref{thm:Pic1} of Appendix \ref{appendix}.

We can also get restrained Pic-categories from chain complexes.   A \define{chain complex of abelian groups} is a sequence of abelian groups $C_0, C_1, \dots, $ together with homomorphisms
\[     C_0 \xleftarrow{\partial_1} C_1 \xleftarrow{\partial_2} C_2 \xleftarrow{\partial_3} \cdots\]
such that $\partial_n \circ \partial_{n+1} = 0$.   A \define{2-term} chain complex of abelian groups is one where $C_n = 0$ except for $C_0$ and $C_1$.  Thus a 2-term chain complex of abelian groups is just an elaborate way of thinking about two abelian groups and a homomorphism between them.  However, this way of thinking is useful because it paves the way for generalizations.

Given a 2-term chain complex $C$ of abelian groups, say $C_0 \xleftarrow{\;\partial\;} C_1 $, we can construct a category $\G_C$ where:
\begin{itemize}
\item  objects are elements $g \in C_0$,
\item  a morphism $h \maps g \to g'$ is an element $h \in C_1$ with $dh = g' - g$,
\item to compose morphisms $h \maps g \to g'$ and $h' \maps g' \to g''$ we add them,
obtaining $h' + h \maps g \to g''$.
\end{itemize}
We can make this category $\G_C$ into a Gr-category by using addition in $C_0$ as the tensor product of objects and addition in $C_1$ as the tensor product of morphisms.  Since addition in these abelian groups is commutative we can make $\G_C$ into a symmetric monoidal category where the symmetry $S_{g,g'} \colon g + g' \to g' + g$ is the identity.  Then $\G_C$ becomes a restrained Pic-category.

In a 1982 paper \cite{H82} Ho\`ang Xu\^an S\'inh showed that in a certain sense \emph{all} restrained Pic-categories arise from this simple construction.  More precisely, she proved that every restrained Pic-category is equivalent to one of the form $\G_C$ for some 2-term chain complex $C$ of abelian groups.   We state this result in Theorem \ref{thm:Pic2} of Appendix \ref{appendix}.  As a consequence, the Pic-category of a commutative ring, which we defined using invertible $R$-modules and isomorphisms between these, can also be expressed in terms of a 2-term chain complex.

\section{Gr-categories and topology}
\label{sec:topology}

For more examples of Gr-categories it pays to exploit the connection between Gr-categories and topology.  Our brief discussion of group cohomology touched on this theme, but did not do justice to either the history or the mathematics.  

One of the dreams of topology is to classify topological spaces.   Two such spaces $X$ and $Y$ are homeomorphic if there are maps going back and forth, say $f \maps X \to Y$ and $g \maps Y \to X$, that are inverses:
\[        gf = 1_X \textrm{\; and \; }  fg = 1_Y. \]    
But it has long been known that classifying spaces up to homeomorphism is an absolutely unattainable goal.  So, various lesser but still herculean tasks have been proposed as substitutes.  

For example, instead of demanding that $f$ and $g$ are inverses `on the nose', we can merely ask for them to be inverses up to
homotopy.   Given two maps $h, h' \maps A \to B$, a \define{homotopy} between them is a continuous 1-parameter family of maps interpolating between them: that is, map $H \maps [0,1] \times A \to B$ with 
\[    \begin{array}{ccl}
H(0,a) &=& h(a),  \\
H(1,a) &=& h'(a) .
\end{array}
\]
If there is a homotopy from $h$ to $h'$ we write $h \simeq h'$.   We then say two spaces $X$ and $Y$ are \define{homotopy equivalent} if there are maps $f \maps X \to Y$ and $g \maps Y \to X$ such that 
\[        gf \simeq 1_X \textrm{\; and \; }  fg \simeq 1_Y. \]  

Classifying spaces up to homotopy equivalence is still impossibly complicated unless we require that these spaces are locally nice in some sense.  Manifolds, being locally homeomorphic to $\R^n$, are the paradigm of what we might mean by locally nice.  But we do not need to go this far.   A `CW complex' is a space built by starting with a discrete set of points, or 0-balls, and iteratively gluing on 1-balls, 2-balls, and so on, by attaching their boundaries to the space built so far.   Topologists have adopted the goal of classifying CW complexes up to homotopy equivalence as a kind of holy grail of the subject.

We can simplify this quest in a couple of ways.  First, since every CW complex is a disjoint union of connected CW complexes, we can focus on connected one.  Second, it is very useful to equip a connected CW complex $X$ with a chosen point, called a `basepoint' and denoted $\ast \in X$.  Spaces with basepoint are called \define{pointed} spaces.  When we working with these we demand that all maps send the basepoint of one space to that of another, and that homotopies do this for all values of the parameter $t$.

Now, instead of trying to classify all connected pointed CW complexes in one blow, which is still beyond our powers, it is more manageable to start with ones whose interesting features are concentrated in low dimensions and work our way up.   For this we should recall the idea of homotopy groups.   Suppose $X$ is a pointed CW complex, and let $S^n$ be the $n$-sphere with an arbitrarily chosen basepoint.  Then we define the \define{$n$th homotopy group} $\pi_n(X)$ to be the set of homotopy classes of maps from $S^n$ to $X$---where again, we demand that all maps and homotopies are compatible with the chosen basepoints.   We can think of $\pi_n(X)$ as the set of holes in $X$ that can be caught with an $n$-dimensional lasso.   Despite its name $\pi_0(X)$ is a mere set, and a 1-element set when $X$ is connected, so this gives no information at all.  But $\pi_n(X)$ is a group for $n \ge 1$, and an abelian group for $n \ge 2$.  

We say a CW complex $X$ is an \define{$n$-type} if, regardless of how we choose a basepoint for it, $\pi_k(X)$ is trivial for $k > n$.  Thus, intuitively, an $n$-type is a locally nice space whose interesting features live only in dimensions $\le n$, at least as viewed through the eye of homotopy groups.  

A remarkable fact, discovered by Eilenberg and Mac Lane around 1954, is that connected pointed 1-types are classified by groups \cite{EM54,EM55}.   On the one hand, given a connected pointed 1-type $X$ we get a group $\pi_1(X)$.   More surprisingly, two connected pointed 1-types have isomorphic groups $\pi_1$ if and only if they are homotopy equivalent.   On the other hand, given any group $G$, there is a concrete procedure for building a connected pointed 1-type $X$ with $\pi_1(X) = G$.   This space $X$ is now called the \define{Eilenberg--Mac Lane space} $K(G,1)$, and it is quite easy to describe.

Start with a point $\ast$ to serve as the basepoint.  For each element $g \in G$, take an interval and glue both of its ends to this point:
\[
\xy 0;/r3pc/:
   (0,0)*+{\ast}="x";
    (2,0)*+{\ast}="z";
 {\ar^{g} "x";"z"};
\endxy
\]
Here we draw the two ends of the interval separately for convenience, but they are really the same point $\ast$. Then, for any pair of elements $g,h \in G$, take a triangle and glue its edges to intervals for $g, h$ and $gh$:
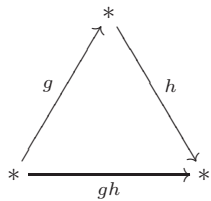
\[
\xy 0;/r3pc/:
   (0,0)*+{\ast}="x";
    (1,1.7)*+{\ast}="y";
    (2,0)*+{\ast}="z";
 {\ar^{g} "x";"y"};
 {\ar^{h} "y";"z"};
 {\ar_{gh} "x";"z"};
\endxy
\]
Next, for any triple of elements $g,h,k \in G$, glue a tetrahedron of this sort to the already present triangles:
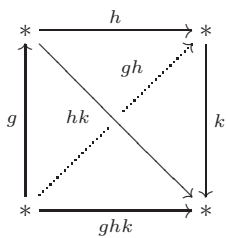
\[
{\xy 
(12,-12)*+{\ast}="1"; 
(12,12)*+{\ast}="2"; 
(-12,12)*+{\ast}="3"; 
(-12,-12)*+{\ast}="4"; 
(0,0)*+{}="0"; 
{\ar^-{k} "2";"1"}; 
{\ar^-{h} "3";"2"}; 
{\ar^-{g} "4";"3"};
{\ar_-{ghk} "4";"1"};
{\ar@{.>}^(0.4){\phantom{ab}gh} "0";"2"};
{\ar@{.} "4";"0"};
{\ar_(0.4){hk} "3";"1"};
\endxy}
\]
And so on: glue on an $n$-dimensional simplex for each $n$-tuple
of elements of $G$.   The result is the Eilenberg--Mac Lane space
$K(G,1)$.

And now we reach the main point: just as connected pointed 1-types are classified by groups, connected pointed 2-types are classified by Gr-categories!  

This realization did not come out of the blue---far from it.  A crucial first step was J.\ H.\ C.\ Whitehead's concept of crossed module, formulated around 1946 without the aid of category theory \cite{W46,W49}.  In 1950, Mac Lane and Whitehead \cite{MW} proved that a crossed module captures all the information in a connected pointed 2-type.   As mentioned in our discussion of strict Gr-categories, it seems Verdier knew by 1965 that a crossed module is another way of thinking about a strict Gr-category (though the latter term did not yet exist).  A proof was first published by Brown and Spencer \cite{BS76} in 1976.  However, Grothendieck was familiar with some of these ideas before then, which must be one reason he proposed that Ho\`ang Xu\^an S\'inh write her thesis on Gr-categories. 

Indeed, Grothendieck already had ideas for vast generalizations, which became more important in his later work, and especially his monumental text
\textsl{Pursuing Stacks} \cite{G84}.  He conjectured that $n$-types should be classified up to homotopy equivalence by algebraic structures called `$n$-groupoids'.  Very roughly speaking, an $n$-groupoid has objects, morphisms between objects, 2-morphisms between 1-morphisms, and so on up to level $n$, all invertible up to higher morphisms.   Filling in the details and proving Grothendieck's conjecture---usually called the `homotopy hypothesis'---continues to be a challenge \cite{HL19}.

An $n$-groupoid with one object is called an `$n$-group'.   A 1-group is simply a group, while a 2-group is a Gr-category.  As a spinoff of the homotopy hypothesis, connected pointed $n$-types should be classified up to homotopy equivalence by $n$-groups.   As we have seen, the case $n = 1$ was handled by Eilenberg and Mac Lane.  While the case $n = 2$ was tackled by the authors listed above, it was still not completely clarified when Ho\`ang Xu\^an S\'inh wrote her thesis.  

In fact, it is possible to mimic Eilenberg and Mac Lane's construction and build a connected pointed 2-type starting from any Gr-category $\G$, much as they built a connected pointed 1-type starting from a group.  For this it is convenient to start by replacing $\G$ with an equivalent strict Gr-category $\G$.   So, let us assume $\G$ is strict.  

Start with a point $\ast$.  For each object $g \in \G$, take an interval and glue both of its ends to this point:
\[
\xy 0;/r3pc/:
   (0,0)*+{\ast}="x";
    (2,0)*+{\ast}="z";
 {\ar^{g} "x";"z"};
\endxy
\]
Then, for any triple of objects $g,h,k \in \G$ and any morphism $a \maps g \otimes h \to k$, take a triangle and glue its edges to the intervals for $g, h$ and $k$:
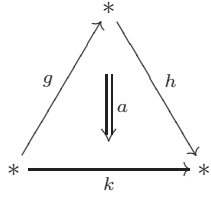
\[
\xy 0;/r3pc/:
   (0,0)*+{\ast}="x";
    (1,1.7)*+{\ast}="y";
    (2,0)*+{\ast}="z";
{\ar@2{->}^-{a} (1,1)*{};(1,0.3)*{}};
 {\ar^{g} "x";"y"};
 {\ar^{h} "y";"z"};
 {\ar_{k} "x";"z"};
\endxy
\]
Next, for any tetrahedron as shown below, glue this tetrahedron to the already present triangles if this tetrahedron commutes:
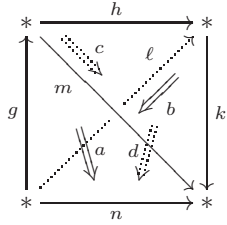
\[
{\xy 
(12,-12)*+{\ast}="1"; 
(12,12)*+{\ast}="2"; 
(-12,12)*+{\ast}="3"; 
(-12,-12)*+{\ast}="4"; 
(0,0)*+{}="0"; 
{\ar^-{k} "2";"1"}; 
{\ar^-{h} "3";"2"}; 
{\ar^-{g} "4";"3"};
{\ar_-{n} "4";"1"};
{\ar@{.>}^(0.5){\tiny{\ell}} "0";"2"};
{\ar@{.} "4";"0"};
{\ar_(0.3){m} "3";"1"};
{\ar@2{->}^-{a} (-5,-2)*{};(-3,-9)*{}};
{\ar@2{->}^-{b} (8,5)*{};(3,0)*{}};
{\ar@2{.>}^(0.6){c} (-7,10)*{};(-2,5)*{}};
{\ar@2{.>}_-{d} (5,-2)*{};(3,-9)*{}};
\endxy}
\]
\noindent
By saying that this tetrahedron `commutes', we mean that the composite
of the front two sides equals the composite of the back two sides, i.e.\ the following two morphisms in $\G$ are equal:
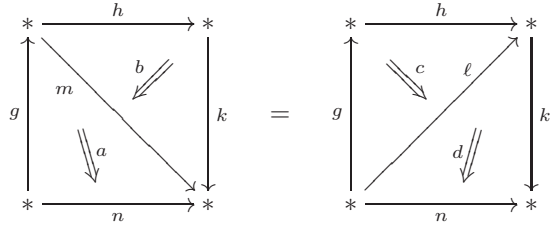
\[
{\xy 
(12,-12)*+{\ast}="1"; 
(12,12)*+{\ast}="2"; 
(-12,12)*+{\ast}="3"; 
(-12,-12)*+{\ast}="4"; 
{\ar^-{k} "2";"1"}; 
{\ar^-{h} "3";"2"}; 
{\ar^-{g} "4";"3"};
{\ar_-{n} "4";"1"};
{\ar_(0.3){m} "3";"1"};
{\ar@2{->}^-{a} (-5,-2)*{};(-3,-9)*{}};
{\ar@2{->}_-{b} (7,7)*{};(2,2)*{}};
\endxy}
\quad = \quad
{\xy 
(12,-12)*+{\ast}="1"; 
(12,12)*+{\ast}="2"; 
(-12,12)*+{\ast}="3"; 
(-12,-12)*+{\ast}="4"; 
{\ar^-{k} "2";"1"}; 
{\ar^-{h} "3";"2"}; 
{\ar^-{g} "4";"3"};
{\ar_-{n} "4";"1"};
{\ar_{} "4";"2"}; 
(3.5,6)*+{\scriptstyle{\ell}};
{\ar@2{->}^(0.5){c} (-7,7)*{};(-2,2)*{}};
{\ar@2{->}_-{d} (5,-2)*{};(3,-9)*{}};
\endxy}
\]
Next  glue on a 4-simplex whenever all its tetrahedral faces are present, and so on in all higher dimensions. This gives a space $B\G$ called the \define{classifying space} of $\G$.  This is the connected pointed 2-type corresponding to the Gr-category $\G$.  For more details see Duskin \cite{D01}, Bullejos and Cegarra \cite{BC03}, and Noori \cite{N08}.

Remember, a Gr-category $\G$ is characterized up to equivalence by four pieces of data:
\begin{itemize}
\item a group $G$,
\item an abelian group $A$,
\item an action $\rho$ of $G$ on $A$,
\item the S\'inh invariant $[a] \in H^3(G,A)$.
\end{itemize}
These are all visible from the classifying space $B\G$.  This space has $\pi_1(B\G) = G$ and $\pi_2(B\G) = A$.  For any pointed space the group $\pi_1$ has an action on the abelian group $\pi_2$, and for the space $B\G$ this action is $\rho$.  More subtly, every space has a `Postnikov invariant', an element of $H^3(\pi_1,\pi_2)$, and for $B\G$ this is the S\'inh invariant $[a] \in H^3(G,A)$.

We have seen how to turn a Gr-category into a space, which happens to be a connected pointed 2-type.  But the relation between Gr-categories and topology runs deeper than that.  There is also a way to turn a pointed space into a Gr-category!  

Suppose $X$ is a pointed space.  First, remember that we can get a group from $X$, namely $\pi_1(X)$, usually called the \define{fundamental group} of $X$.   Elements of $\pi_1(X)$ are homotopy classes of loops in $X$: that is, maps $g \maps S^1 \to X$.   To multiply two loops we simply combine them into a single loop, starting the second where the first ended.    

In Example 2.3.3 of her thesis \cite{H73}, Ho\`ang Xu\^an S\'inh tersely described how to enhance this procedure to get a Gr-category from $X$.  This Gr-category is now called the \define{fundamental 2-group} of $X$, and denoted $\Pi_2(X)$.   An object of $\Pi_2(X)$ is a map $g \maps S^1 \to X$, and a morphism $f \maps g \to g'$ is a homotopy class of homotopies from $g$ to $g'$.  The details are carefully worked out in various later sources, e.g.\ the work of Hardie, Kamps and Kieboom \cite{HKK00,HKK01}. 

There are other ways to obtain the fundamental 2-group of a connected pointed space $X$.  For example:
\begin{itemize}
\item let $G = \pi_1(X)$,
\item let $A = \pi_2(X)$,
\item let $\rho$ be the action of $\pi_1(X)$ on $\pi_2(X)$,
\item let $[a] \in H^3(\pi_1(X),\pi_2(X))$ be the Postnikov invariant of $X$.
\end{itemize}
Then, we can take $\Pi_2(X)$ to be the Gr-category associated to this quadruple $(G,A,\rho, [a])$.  

As a result of all this, not only can we turn any Gr-category $\G$ into a 
connected pointed 2-type $B\G$, we can also turn any connected pointed 2-type $X$ into a Gr-category $\Pi_2(X)$.   Furthermore, with work 
one can show that:
\begin{itemize}
\item Any Gr-category $\G$ is equivalent to $\Pi_2(B\G)$.
\item Any connected pointed 2-type $X$ is homotopy equivalent to 
$B(\Pi_2(X))$.
\end{itemize}
As a corollary, connected pointed 2-types are classified by Gr-categories.  That is, two connected pointed 2-types are homotopy equivalent if and only if their fundamental 2-groups are equivalent as Gr-categories.

But one can equally well turn the logic around and say Gr-categories are classified by connected pointed 2-types.  Two Gr-categories are equivalent if and only if their classifying spaces are homotopy equivalent.   So, the wall separating algebra and topology has been completely knocked down, at least in this limited realm!  This is a nice piece of evidence for Grothendieck's homotopy hypothesis.

\section{Gr-categories in physics}
\label{sec:physics}

While our story has been focused on Ho\`ang Xu\^an S\'inh's thesis and its immediate context, it would be a disservice not to mention that Gr-categories---now often called 2-groups \cite{BL04}---later took on a new life in theoretical and mathematical physics.   The reason is that `gauge theory', the spectacularly successful approach to physics based on groups, can be generalized to `higher gauge theory' using 2-groups \cite{BH11,BS06}.  Just as gauge theory based on groups describes how point particles change state as they trace out 1-dimensional paths in spacetime, higher gauge theory based on 2-groups describes both particles and strings, the latter of which trace out 2-dimensional `worldsheets' in spacetime.  Furthermore, there is no reason to stop with 2-groups \cite{NSS23,S13}.    While no theories of physics based on higher gauge theory have received any experimental confirmation, the mathematics behind these theories is magnificent, and there are hopes that someday physicists will synthesize phases of matter described by higher gauge theory \cite{BCHK22,BHKR21}.  This would be a wonderful realization of Ho\`ang Xu\^an S\'inh's vision.

\appendix
\section{Equivalence for Gr-categories and Pic-categories}
\label{appendix}

As we have seen, the concept of `equivalence' is crucial for many of 
Ho\`ang Xu\^an S\'inh's results on Gr-categories and Pic-categories.  This concept is subtler than isomorphism.  We can see this already for categories.  Two categories $\mathsf{C}$ and $\mathsf{D}$ are isomorphic if there are functors going back and forth, $F \maps \mathsf{C} \to \mathsf{D}$ and $G \maps \mathsf{D} \to \mathsf{C}$, that are inverses: $GF = 1_{\mathsf{C}}$ and $FG = 1_\mathsf{D}$.    But experience has taught us than demanding \emph{equations} between functors is too harsh: instead it suffices to have \emph{natural isomorphisms} between them.  So, we say $\mathsf{C}$ and $\mathsf{D}$ are \define{equivalent} if there are functors $F \maps \mathsf{C} \to \mathsf{D}$ and $G \maps \mathsf{D} \to \mathsf{C}$ for which there exist natural isomorphisms
\[   \alpha \maps GF \stackrel{\sim}{\Longrightarrow} 1_{\mathsf{C}}, \quad \beta \maps FG \stackrel{\sim}{\Longrightarrow} 1_\mathsf{D}. \]
(We use a single arrow for functors and a double arrow for natural transformations.)  For Gr-categories and Pic-categories, we need to enhance this concept of equivalence to take into account the extra structure that these categories carry.
 
As we have seen, Gr-categories are simply monoidal categories with extra properties.  Similarly Pic-categories are symmetric monoidal categories with extra properties.  So, it is enough to define notions of equivalence for monoidal and symmetric monoidal categories.  We begin by defining monoidal functors.  These do not need to preserve the tensor product and unit object strictly, but rather only up to natural isomorphism.  However, these natural isomorphisms need to obey some laws of their own!

\begin{defn} 
A functor $F \maps \M \to \N$ between monoidal categories is
\define{monoidal} if it is equipped with: 
\begin{itemize}
\item
a natural isomorphism
$\Phi_{x,y} \maps F(x) \otimes F(y) \xrightarrow{\;\sim\;} F(x \otimes y)$ and
\item
an isomorphism $\phi \maps I_\N \xrightarrow{\;\sim\;} F(I_\M)$
\end{itemize}
such that:
\begin{itemize}
\item
the following diagram commutes for any objects $x,y,z \in \M$:
\[
\xymatrix{
(F(x) \otimes F(y)) \otimes F(z)
\ar[rr]^{\Phi_{x,y}\, \otimes 1_{F(z)}} 
\ar[d]_{a_{F(x),F(y),F(z)}} & &
F(x \otimes y) \otimes F(z)
\ar[rr]^{\Phi_{x\otimes y,z}} & &
F((x \otimes y)\otimes z)
\ar[d]^{F(a_{x,y,z})}  \\
F(x) \otimes (F(y) \otimes F(z))
\ar[rr]_{1_{F(x)} \otimes \Phi_{y,z}} & &
F(x) \otimes F(y \otimes z)
\ar[rr]_{\Phi_{x,y\otimes z}} & &
F(x \otimes (y \otimes z)),
}
\]
\item
the following diagrams commute for any object $x \in \M$:
\[
\xymatrix{
I_\N \otimes F(x) \ar[r]^{\ell_{F(x)}}
\ar[d]_{\phi \otimes 1_{F(x)}} &
F(x) 
&&
F(x) \otimes I_\N \ar[r]^{r_{F(x)}}
\ar[d]_{1_{F(x)} \otimes \phi} &
F(x) 
\\
F(I_\M) \otimes F(x) \ar[r]_{\Phi_{I_\M,x}} &
F(I_\M \otimes x) \ar[u]_{F(\ell_x)}
&&
F(x) \otimes F(I_\M) \ar[r]_{\Phi_{x,I_\M}} &
F(x \otimes I_\M). \ar[u]_{F(r_x)}
}
\]
\end{itemize}
\end{defn}

We similarly have a concept of `symmetric monoidal functor'.  Here the natural isomorphism $\Phi$ must get along with the symmetry:

\begin{defn} 
A functor $F \maps \M \to \N$ between symmetric monoidal categories is
\define{symmetric monoidal} if it is monoidal and it makes the following diagram commute for all $x,y \in \M$:
\[  
\xymatrix{
F(x) \otimes F(y) \ar[rr]^{S_{F(x),F(y)}}
\ar[d]_{\Phi_{x,y}} &&
F(y) \otimes F(x)
\ar[d]^{\Phi_{y,x}} \\
F(x \otimes y) \ar[rr]_{F(S_{x,y})}  &&
F(y \otimes x)
}
\]
\end{defn}  

Next we need monoidal \emph{natural transformations}.  Recall that a monoidal functor $F \maps \M \to \N$ is really a triple $(F, \Phi, \phi)$ consisting of a functor, a natural isomorphism $\Phi$, and an isomorphism $\phi$.  A `monoidal natural transformation' must get  along with these extra isomorphisms: 

\begin{defn} 
Suppose that $(F,\Phi,\phi)$ and $(G,\Gamma,\gamma)$ are monoidal
functors from the monoidal category $\M$ to the monoidal category $\N$. 
Then a natural transformation $\alpha \maps F \To G$ is \define{monoidal} if 
the diagrams
\[ 
\xymatrix{
F(x) \otimes F(y) 
\ar[rr]^{\alpha_x \otimes \alpha_y}
\ar[d]_{\Phi_{x,y}} &&
G(x) \otimes G(y) \ar[d]^{\Gamma_{x,y}}  
\\
F(x \otimes y)  \ar[rr]_{\alpha_{x\otimes y}} &&
G(x \otimes y)
}
\]
and
\[
\xy 0;/r.3pc/:
(0,0)*+{I_\N}="t";
(-10,-10)*+{F(I_\M)}="l";
(10,-10)*+{G(I_\M)}="r";
{\ar_\phi "t";"l"};
{\ar^\gamma "t";"r"};
{\ar_{\alpha_{I_\M}} "l";"r"};
\endxy    
\]
commute.
\end{defn}
\noindent There are no extra conditions required of \define{symmetric monoidal} natural transformations: they are simply monoidal natural transformations between symmetric monoidal functors.

We can compose two monoidal natural transformations and get another monoidal natural transformation.  Also, any monoidal functor $F$ has an identity monoidal natural transformation $1_F \maps F \To F$.    This allows us to make the following definition:

\begin{defn}
A monoidal natural transformation $\alpha \maps F \To G$ is a
\define{monoidal natural isomorphism} if there is a monoidal natural transformation $\beta \maps G \To F$ that is an inverse to $\alpha$: 
\[       \beta \alpha = 1_F   \textrm{ \; and \; } \alpha \beta = 1_G .\]
\end{defn}

Now at last we are ready to define equivalence for Gr-categories and Pic-categories, so we can state some of Ho\`ang Xu\^an S\'inh's major results in a precise way.  For this we need the fact that we can compose two monoidal functors and get a monoidal functor, and similarly in the symmetric monoidal case.

\begin{defn}  
Two monoidal categories $\M$ and $\N$ are \define{equivalent} if there exist
monoidal functors $F \maps \M \to \N$ and $G \maps \N \to \M$ such that 
there exist monoidal natural isomorphisms $\alpha \maps GF \To 1_\N$ and 
$\beta \maps FG \To 1_\N$.   Two Gr-categories $\G$ and $\H$ are equivalent if they are equivalent as monoidal categories.
\end{defn}

\begin{defn}  
\label{defn:Pic_equivalence}
Two symmetric monoidal categories $\M$ and $\N$ are \define{equivalent} if there exist symmetric monoidal functors $F \maps \M \to \N$ and $G \maps \N \to \M$ such that there exist symmetric monoidal natural isomorphisms $\alpha \maps GF \To 1_\N$ and $\beta \maps FG \To 1_\N$.  Two Pic-categories $\G$ and $\H$ are equivalent if they are equivalent as symmetric monoidal categories.
\end{defn}

With these concepts in hand, we can state some results on the classification of Gr-categories and Pic-categories.   We have already mentioned the following pair of results.  The former is Prop.\ 3.22 in Ho\`ang Xu\^an S\'inh's thesis \cite{H73}, while the latter is Prop.\ 5.3 in her 1978 paper on strict Gr-categories \cite{H78}.

\begin{thm}
\label{thm:Gr1}
Any Gr-category is equivalent to a skeletal Gr-category, i.e.\ one for which
isomorphic objects are necessarily equal.
\end{thm}

\begin{thm}
\label{thm:Gr3}
Any Gr-category is equivalent to a strict Gr-category, i.e.\ one for which the associator, left unitor and right unitor are all identity natural transformations.
\end{thm}

We have also seen how any skeletal Gr-category $\G$ gives a 4-tuple $(G,A,\rho,\alpha)$.  The following result gives a classification of skeletal Gr-categories in these terms.  It is a special case of Proposition 3.47 of Ho\`ang Xu\^an S\'inh's thesis \cite{H73}:

\begin{thm}
\label{thm:Gr2}
Suppose $\G$ and $\G'$ are skeletal Gr-categories giving 4-tuples $(G,A,\rho,\alpha)$ and $(G',A',\rho',\alpha')$, respectively.  Then $\G$ and $\G'$ are equivalent if and only if there exist isomorphisms $\phi \maps G \to G'$ and $\psi \maps A \to A'$ such that:
\begin{itemize}
\item
the actions $\rho$ and $\rho'$ are related as follows:
\[   \rho'(\phi(g))(\psi(a)) = \psi(\rho(g)(a)) \]
for all $g \in G$ and $a \in A$,
\item the cohomology classes $\alpha$ and $\alpha'$ are related as follows:
\[   \alpha'(\phi(g_1),\phi(g_2),\phi(g_3)) - \psi(\alpha(g_1,g_2,g_3)) = df \]
for some $f \maps G'^2 \to A$.
\end{itemize}
\end{thm}
\noindent
She stated a more general result that does not require $\G$ and $\G'$ to be skeletal, but since every Gr-category is equivalent to a skeletal one, a classification of skeletal Gr-categories up to equivalence can serve as a classification of all Gr-categories.   For a deeper treatment of the relation between Gr-categories and group cohomology, see Section 6 of Joyal and Street's unpublished paper \cite{JS86}.

Interestingly, the unitors do not provide any extra information required for the classification of Gr-categories up to equivalence.  The reason is that given any Gr-category $\G$, we can change the tensor product by setting $I \otimes g = g = g \otimes I$ for every $g \in \G$, change the left and right unitors to identity morphisms, and adjust the associators to obtain a new Gr-category $\G'$ that equivalent to $\G$.  For a self-contained proof see \cite[Prop.\ 39]{BL04}.  

The following classification of restrained Pic-categories is Corollary 3.61 in Ho\`ang Xu\^an S\'inh's thesis \cite{H73}. She derived it as a consequence of a subtler classification of \emph{all} Pic-categories, her Proposition 3.59.

\begin{thm}
\label{thm:Pic1}
Suppose $\G$ and $\H$ are restrained Pic-categories.  Then $\G$ and $\H$ are equivalent if and only if both these conditions hold:
\begin{itemize}
\item
The abelian group of isomorphism classes of objects in $\G$ is isomorphic to the
abelian group of isomorphism classes of objects in $\H$.
\item
The abelian group of automorphisms of $I_\G$ is isomorphic to the abelian
group of automorphisms of $I_\H$.
\end{itemize}
\end{thm}

Later, in her 1982 paper ``Cat\'egories de Picard restreintes'' \cite{H82}, Ho\`ang Xu\^an S\'inh showed that every restrained Gr-category arises from a 2-term chain complex of abelian groups in the manner described at the end of Section \ref{sec:Pic}:

\begin{thm}
\label{thm:Pic2}
Suppose $\G$ is a restrained Pic-category.  Then there is some 2-term chain complex $C$ of abelian groups such that $\G$ is equivalent, as a Pic-category, to $\G_C$.
\end{thm}

\vskip 1em

\subsection*{Acknowledgements}   

I thank Ho\`ang Xu\^an S\'inh for answering my questions about her thesis, Lisa
Raphals for translating my questions into French, and David Egolf and David Michael Roberts for pointing out mistakes and ways in which this paper could be improved.

\end{document}